\documentclass[12pt]{amsart}
\usepackage[]{amscd,amsfonts,amsmath,amsxtra,amssymb}
\usepackage{graphics}
\usepackage[all]{xy}
\usepackage{hyperref}
\usepackage{eucal}
\newtheorem{sub}{}[section]
\newtheorem{subsub}{}[sub]
\def\ov#1{\overline{#1}}

\def\HHom{\mathop{\mathcal Hom}\nolimits}

\def\Ext{\mathop{\rm Ext}\nolimits}

\def\Pic{\mathop{\rm Pic}\nolimits}

\def\deg{\mathop{\rm deg}\nolimits}

\def\spec{\mathop{\rm spec}\nolimits}
\def\lra{\longrightarrow}

\def\sigg{\mathop{\hbox{$\displaystyle\sum$}}\limits}

\def\hfl#1#2{\smash{\mathop{\ \hbox to 12mm{\rightarrowfill}}
\limits^{\scriptstyle#1}_{\scriptstyle#2} \ }}
\def\hflb#1#2{\smash{\mathop{\hbox to 12mm{\leftarrowfill}}
\limits^{\scriptstyle#1}_{\scriptstyle#2}}}

\def\m#1{{\hbox{$#1$}}}

\def\ot{\otimes}
\def\og{\leavevmode\raise.3ex\hbox{$\scriptscriptstyle\langle\!\langle$}}
\def\fg{\leavevmode\raise.3ex\hbox{$\scriptscriptstyle\,\rangle\!\rangle$}}
\def\span#1{\langle#1\rangle}

\def\nsp{\lbrace 0\rbrace}

\def\dsp{\displaystyle}
\def\Ssect#1#2{\pagebreak[3]\begin{sub}\label{#2}{\sc\small\small
#1}\rm\medskip}

\def\sepsec{\vskip 1.4cm}
\def\sepsub{\vskip 0.7cm}

\def\sepprop{\vskip 0.5cm}

\def\xmat#1{\[\xymatrix{#1}\]}
\def\flinc{\ar@{^{(}->}}
\def\fleq{\ar@{=}}
\def\flon{\ar@{->>}}
\def\fmaps{\ar@{|-{>}}}
\def\Nligne{\hfil\break}

\def\BS#1{{\boldsymbol{#1}}}

\newcommand{\C}{{\mathbb C}}

\renewcommand{\P}{{\mathbb P}}

\newcommand{\ka}{{\mathcal A}}
\newcommand{\kb}{{\mathcal B}}
\newcommand{\kc}{{\mathcal C}}
\newcommand{\kd}{{\mathcal D}}
\newcommand{\ke}{{\mathcal E}}
\newcommand{\kf}{{\mathcal F}}

\newcommand{\ki}{{\mathcal I}}

\newcommand{\ko}{{\mathcal O}}

\newcommand{\kx}{{\mathcal X}}

\newcommand{\kz}{{\mathcal Z}}
\begin{document}

\def\refname{References}
\def\contentsname{Summary}
\def\proofname{Proof}
\def\abstractname{Resume}

\author{Jean--Marc Dr\'{e}zet}
\address{
Institut de Math\'ematiques de Jussieu,
Case 247,
4 place Jussieu,
F-75252 Paris, France}
\email{jean-marc.drezet@imj-prg.fr}
\title[{Reducible deformations and smoothing}] {Reducible deformations and
smoothing of primitive multiple curves}

\begin{abstract}
A {\em primitive multiple curve} is a Cohen-Macaulay irreducible projective
curve $Y$ that can be locally embedded in a smooth surface, and such that
$C=Y_{red}$ is smooth. In this case, $L=\ki_C/\ki_C^2$ is a line bundle on $C$.

This paper continues the study of deformations of $Y$ to curves with
smooth irreducible components, when the number of components is maximal (it is
then the multiplicity $n$ of $Y$). If a primitive double curve $Y$ can
be deformed to reduced curves with smooth components intersecting 
transversally, then $h^0(L^{-1})\not=0$. We prove that conversely, if $L$ is 
the ideal sheaf of a divisor with no multiple points, then $Y$ can
be deformed to reduced curves with smooth components intersecting 
transversally. We give also some properties of reducible
deformations in the case of multiplicity $n>2$.
\end{abstract}

\maketitle
\tableofcontents

\vskip 1cm

\section{Introduction}\label{intro}

A {\em primitive multiple curve} is an algebraic variety $Y$ over $\C$ which is
Cohen-Macaulay, such that the induced reduced variety \m{C=Y_{red}} is a
smooth projective irreducible curve, and that every closed point of $Y$ has a
neighborhood that can be embedded in a smooth surface. These curves have been
defined and studied by C.~B\u anic\u a and O.~Forster in \cite{ba_fo}. The
simplest examples are infinitesimal neighborhoods of projective smooth curves
embedded in a smooth surface (but most primitive multiple curves cannot be
globally embedded in smooth surfaces). Primitive multiple curves of 
multiplicity 2 (called {\em ribbons}) have been parametrized in \cite{ba_ei}. 
Primitive multiple curves of any multiplicity and the coherent sheaves on them 
have been studied in \cite{dr1}, \cite{dr2}, \cite{dr4} and \cite{dr5}.

Let $Y$ be a primitive multiple curve with associated reduced curve $C$, and
suppose that \m{Y\not=C}. Let \m{\ki_C} be the ideal sheaf of $C$ in $Y$. The
{\em multiplicity} of $Y$ is the smallest integer $n$ such that \m{\ki_C^n=0}.
We have then a filtration
\[C=C_1\subset C_2\subset\cdots\subset C_{n}=Y \]
where \m{C_i} is the subscheme corresponding to the ideal sheaf \m{\ki_C^i}
and is a primitive multiple curve of multiplicity $i$. The sheaf \ \m{L=
\ki_C/\ki_C^2} \ is a line bundle on $C$, called the {\em line bundle on $C$
associated to $Y$}.

\sepsub

\Ssect{Deformations to reduced reducible curves}{def1}

Deformations of primitive multiple curves \m{Y=C_n} of any multiplicity
\m{n\geq 2} to reduced curves having multiple components which are smooth,
intersecting transversally, have been studied in \cite{dr7}: $n$ is the maximal
number of components of such deformations of $Y$, and in this case we say that
the deformation is {\em maximal}, and the number of intersection points of two
components is exactly \m{-\deg(L)} (these deformations are called {\em maximal 
reducible deformations}). In \cite{dr7} the case \m{\deg(L)=0} has
been completely treated: a primitive multiple curve of multiplicity $n$ can be
deformed in disjoint unions of $n$ smooth curves if and only if \m{\ki_C} is
isomorphic to the trivial bundle on \m{C_{n-1}}. In \cite{dr6} it has been
proved that this last condition is equivalent to the following: there exists a
flat family of smooth curves \m{\kc\to S}, parametrized by a smooth curve $S$,
\m{s_0\in S} such that \m{\kc_{s_0}=C}, such that $Y$ is isomorphic to the
$n$-th infinitesimal neighborhood of $C$ in $\kc$.

The problem of determining which primitive multiple curves of multiplicity $n$
can be deformed to reduced curves having exactly $n$ components, allowing
intersections of the components, is more difficult. A necessary condition is \ 
\m{h^0(L^*)>0}. 

We treat mainly the case of primitive double curves. We introduce the notion of 
{\em local reducible deformation} (that can be obtained for example from a 
global deformation by restricting it to a suitable infinitesimal neighbourhood 
of $C$). We prove that if \ \m{Y=C_2} \ is such that \ \m{h^0(L^*)>0}, then it 
admits a local reducible deformation (theorem 3.2.6 (i)). We conjecture that 
any local reducible deformation can be extended to a global one. We prove this 
in the following case: there exist distinct points \m{P_1,\ldots,P_p\in C} such 
that \ \m{L=\ko_C(-P_1-\cdots-P_p)} (theorem 3.2.6 (ii)). So in this case $Y$ 
admits a global maximal reducible deformation.

The simplest primitive double curve (with associated smooth curve $C$ and 
associated line bundle $L$ on $C$) is the {\em trivial curve} : to obtain it we 
embed $C$ in $L^*$ (a smooth surface) using the zero section, and take for $Y$ 
the second infinitesimal neigbourhood of $C$ in $L^*$. But there is no similar 
easy way to construct a general primitive double curve. Most of them cannot be 
embedded in a smooth surface: in the case \m{C=\P_1}, D.~Bayer and D.~Eisenbud
have proved in \cite{ba_ei}, theorem 7.1, that the only primitive double curves
 that can be embedded in a smooth surface are the trivial ones, and the conic 
 in $\P_2$. But the non trivial primitive double curves (corresponding to 
 $\P_1$ and $L$) are parametrized by the projective space \m{\P(H^1(L(2))}, so 
 there are many non embeddable primitive double curves when \m{\deg(L)\leq 5}.
So to treat the general case we have to use a more abstract construction of 
primitive double curves. We use the parametrization of primitive double curves: 
it was obtained for the first time in \cite{ba_ei}, and in another way in 
\cite{dr1}, using \v Cech cohomology : primitive double curves are constructed 
by gluing smaller pieces of the form \ \m{U\times\spec(\C[t]/(t^2))}, where $U$ 
is an open subset of $C$. This construction can also be used to construct 
local reducible deformations (at least when there exist distinct points 
\m{P_1,\ldots,P_p\in C} such that \ \m{L=\ko_C(-P_1-\cdots-P_p)}). This is done 
in chapter \ref{Maxred2}.

In chapter \ref{Maxredn}, the properties of reducible deformations of double
curves used in chapter \ref{Maxred2} are extended to the case of multiplicity
\m{n>2}. This could be useful to determine which primitive multiple curves of
multiplicity $n$ can be deformed to reducible curves with $n$ smooth
components, using the parametrization of primitive multiple curves given in
\cite{dr1}.

\end{sub}

\sepsub

\Ssect{Smoothing of primitive multiple curves}{smo1}

Since curves with smooth components intersecting transversally are smoothable
(cf. \ref{arg2}), any primitive multiple curve having a maximal reducible
deformation is smoothable.

The deformations of primitive double (i.e. of multiplicity 2) curves
(also called {\em ribbons}) to smooth projective curves have been studied by
M.~Gonz\'alez in \cite{gon}: he proved that such a curve $Y$, with associated
smooth curve $C$ and associated line bundle $L$ on $C$ is smoothable if \
\m{h^0(L^{-2})\not=0}. Here we prove that $Y$ can be deformed to curves
with 2 components intersecting transversally only if \ \m{h^0(L^{-1})\not=0}. 
So there exist smoothable primitive double curves that cannot be deformed to 
curves with 2 components intersecting transversally.

\end{sub}

\sepsub

\Ssect{Deformations of coherent sheaves}{defcoh}

Another motivation for the study of the deformations of primitive multiple
curves in reducible ones is the understanding of the moduli spaces of
semi-stable coherent sheaves on primitive multiple curves. A coherent sheaf
$\ke$ on \m{Y=C_n} is not in general locally free on some nonempty open subset
of $C$. It has been proved in \cite{dr2} that on some nonempty open subset of
$C$, there exist uniquely determined integers \m{m_i\geq 0}, \m{1\leq i\leq n},
such that $\ke$ is locally isomorphic to a sheaf of the form
\[m_1\ko_C\oplus m_2\ko_{C_2}\oplus\cdots\oplus m_n\ko_{C_n} \ . \]
In some cases, nonempty components of moduli spaces of semi-stable
sheaves on $Y$ contain no generically locally free sheaves (cf. \cite{dr2},
\cite{dr4}, \cite{dr5}). I conjecture that if \m{C_n} has a maximal reducible
deformation \m{\pi:\kc\to S}, then this kind of sheaf can be deformed to
sheaves $\kf$ on the fibers \m{\kc_s} such that the ranks of the restrictions
of $\kf$ to the components of \m{\kc_s} are not the same, and determined by
the integers \m{m_i}. Semi-stable vector bundles on curves with many components
have already been studied in \cite{tei1}, \cite{tei2}, \cite{tei3}.
\end{sub}

\sepsub

I am grateful to the anonymous referee for pointing some incomplete proofs in 
3.2 and for example \ref{exa}.

\sepsub

{\bf Notations : } If $X$ is an algebraic variety and \m{Y\subset X} a 
subvariety, \m{\ki_{Y,X}} (or \m{\ki_Y} if there is no risk of confusion) 
denotes 
the ideal sheaf of $Y$ in $X$.

\newpage

\section{Preliminaries}

\Ssect{Primitive multiple curves}{cmpr}

(cf. \cite{ba_fo}, \cite{ba_ei}, \cite{dr2}, \cite{dr1}, \cite{dr4},
\cite{dr5},
\cite{dr6}, \cite{ei_gr}).

Let $C$ be a smooth connected projective curve. A {\em multiple curve with
support $C$} is a Cohen-Macaulay scheme $Y$ such that \m{Y_{red}=C}.

Let $n$ be the smallest integer such that \m{Y=C^{(n-1)}}, \m{C^{(k-1)}}
being the $k$-th infinitesimal neighborhood of $C$, i.e. \
\m{\ki_{C^{(k-1)}}=\ki_C^{k}} . We have a filtration \ \m{C=C_1\subset
C_2\subset\cdots\subset C_{n}=Y} \ where $C_i$ is the biggest Cohen-Macaulay
subscheme contained in \m{Y\cap C^{(i-1)}}. We call $n$ the {\em multiplicity}
of $Y$.

We say that $Y$ is {\em primitive} if, for every closed point $x$ of $C$,
there exists a smooth surface $S$, containing a neighborhood of $x$ in $Y$ as
a locally closed subvariety. In this case, \m{L=\ki_C/\ki_{C_2}} is a line
bundle on $C$ and we have \ \m{\ki_{C_j}=\ki_C^j},
\m{\ki_{C_{j}}/\ki_{C_{j+1}}=L^j} \ for \m{1\leq j<n}. We call $L$ the line
bundle on $C$ {\em associated} to $Y$. Let \m{P\in C}. Then there exist
elements $y$, $t$ of \m{m_{S,P}} (the maximal ideal of \m{\ko_{S,P}}) whose
images in \m{m_{S,P}/m_{S,P}^2} form a basis, and such that for \m{1\leq i<n}
we have \ \m{\ki_{C_i,P}=(y^{i})} .

The simplest case is when $Y$ is contained in a smooth surface $S$. Suppose
that $Y$ has multiplicity $n$. Let \m{P\in C} and \m{f\in\ko_{S,P}}  a local
equation of $C$. Then we have \ \m{\ki_{C_i,P}=(f^{i})} \ for \m{1<j\leq n},
in particular \m{I_{Y,P}=(f^n)}, and \ \m{L=\ko_C(-C)} .


For any \m{L\in\Pic(C)}, the {\em trivial primitive curve} of multiplicity $n$,
with induced smooth curve $C$ and associated line bundle $L$ on $C$ is the
$n$-th infinitesimal neighborhood of $C$, embedded by the zero section in the
dual bundle $L^*$, seen as a surface.

\end{sub}

\sepsub

\Ssect{Construction of primitive multiple curves}{constpmc}

(cf. \cite{dr1}).

Let \m{(U_i)_{i\in I}} be an affine open cover of $C$. Let $n$ be a
positive integer and \ \m{Z_n=\spec(\C[t]/(t^n))}. For
\m{i,j\in I}, \m{i\not=j}, let
\[\Psi_{ij}:U_{ij}\times Z_n\lra U_{ij}\times Z_n\]
be an automorphism leaving \m{U_{ij}} invariant. Suppose that these
automorphisms verify the cochain condition \
\m{\Psi_{ik}=\Psi_{jk}\circ\Psi_{ij}}. Then by gluing the \m{U_i\times Z_n}
using the \m{\Psi_{ij}} we obtain a primitive multiple curve \m{C_n} of
multiplicity $n$ whose associated reduced curve is $C$. Every primitive
multiple curve can be obtained in this way.

Let \m{U\subset C} be a proper open subset. Suppose that \m{\omega_C} is
trivial on \m{U}. Let \m{x\in\ko_C(U)} such that \m{dx} is a generator of
\m{\omega_C(U)}. Then the automorphisms of \ \m{U\times Z_n} \ leaving $U$
invariant are the same as the automorphisms $\phi$ of the $\C$-algebra
\m{\ko_C(U)[t]/(t^n)} such that for every \m{\gamma\in\ko_C(U)[t]/(t^n)},
\m{\phi(\gamma)_{|C}=\gamma_{|C}}. Such automorphisms are of the form
\m{\phi_{\mu,\nu}}, with \m{\mu,\nu\in\ko_C(U)[t]/(t^{n-1})}, $\nu$ invertible,
defined by: for every \m{\alpha\in\ko_C(U)}
\[\phi_{\mu,\nu}(\alpha) \ = \ \sigg_{i=0}^{n-1}(\mu
t)^i\frac{\partial^i\alpha}{\partial^i x} \ , \]
and
\[\phi_{\mu,\nu}(t) \ = \ \nu t \ . \]
Suppose that \m{\omega_C} is trivial on each \m{U_{ij}} and that we have fixed
a generator \ \m{dx=dx_{ij}\in\omega_C(U_{ij})}. Then we can write \
\m{\Psi_{ij}=\phi_{\mu_{ij},\nu_{ij}}}, with \
\m{\mu_{ij},\nu_{ij}\in\ko_C(U_{ij})[t]/(t^{n-1})}, \m{\nu_{ij}} invertible.
The family \m{(\nu_{ij|C})} is a cocycle defining the line bundle $L$
associated to \m{C_n}. If the ideal sheaf \m{\ki_{C,C_n}} of $C$ in \m{C_n} is
isomorphic to the trivial line bundle on \m{C_{n-1}} then we can assume that
\m{\nu_{ij}=1} for all \m{i,j}.

\sepprop

\begin{subsub}\label{dblcst} The case of double curves -- \rm
Now we suppose that \m{n=2}. In this case we have \m{\mu,\nu\in\ko_C(U)}, $\nu$
is invertible, and \m{\phi_{\mu,\nu}} can also been represented by a matrix
\m{\dsp\begin{pmatrix} 1 & 0\\ \mu\frac{\partial}{\partial x} &
\nu\end{pmatrix}} in such a way that the composition of morphisms is equivalent
to the multiplication of matrices, i.e we have
\[\phi_{\mu',\nu'}\circ\phi_{\mu,\nu} \ = \ \phi_{\mu'',\nu''} \ , \]
with \ \m{\mu''=\mu'+\nu'\mu}, \m{\nu''=\nu\nu'}.

We can see \m{\phi_{\mu_{ij},\nu_{ij}}} as a matrix
\[\begin{pmatrix}1 & 0 \\ \mu_{ij}\frac{\partial}{\partial x_{ij}} & \nu_{ij}
\end{pmatrix} \ . \]
Here \m{(\nu_{ij})} is a cocycle representing $L$ (the line bundle on $C$
associated to \m{C_2}), and \m{(\mu_{ij}\frac{\partial}{\partial x_{ij}})} is a
cocycle representing an element $\eta$ of \ \m{H^1(T_C\ot
L)=\Ext^1(\omega_C,L)}, associated to the canonical exact sequence
\[0\lra L\lra\Omega_{Y\mid C}\lra\omega_C\lra 0 \ . \]
Moreover, \m{C_2} is trivial if and only if \m{\eta=0}.

According to \cite{ba_ei} and \cite{dr1}, \m{\C\eta} defines completely \m{C_2}.
More precisely, we say two primitive double curves \m{C_2}, \m{C'_2}, with the
same induced smooth curve $C$, are {\em isomorphic} if there exists an
isomorphism \m{C_2\simeq C'_2} inducing the identity on $C$. Of course in this
case the associated line bundles on $C$ are also the same. By associating
\m{\C\eta} to \m{C_2}, we define a bijection from the set of isomorphism
classes of non trivial primitive double curves with induced smooth curve $C$
and associated line bundle $L$ and \m{\P(H^1(T_C\ot L))}.

Let $S$ be a smooth curve and \m{P\in C}. Let \m{\pi:\kc\to S} be a flat family 
of projective irreducible smooth curves parametrized by $S$, and suppose that \ 
\m{C=\pi^{-1}(P)}. Let \m{C_2} be the second infinitesimal neighbourhood of $C$ 
in $\kc$, which is a primitive double curve with associated line bundle 
\m{\ko_C}. To this double curve we have associated \m{\eta\in H^1(T_C)}.
Then the image of the Koda\" ira-Spencer map of the deformation $\kc$ of $C$
\[T_PS\lra H^1(T_C)\]
is \m{\C\eta}.
\end{subsub}

\end{sub}

\sepsub

\Ssect{Maximal reducible deformations}{Maxred_dr7}

We recall here some definitions and results of \cite{dr7}.

\begin{subsub}\label{p1} \rm
Let $C$ be a projective irreducible smooth curve, \m{n\geq 2} an integer and
\m{C_n} a primitive multiple curve of multiplicity $n$, with underlying smooth
curve $C$. Let $S$ be a smooth curve, \m{P\in S} and \ \m{\pi:\kc\to S} \
a {\em maximal reducible deformation of \m{C_n}} (cf. \cite{dr7}). This means
that
\begin{enumerate}
\item[(i)] $\kc$ is a reduced algebraic variety with $n$ irreducible
components $\kc_1,\ldots,\kc_n$.
\item[(ii)] We have \ $\pi^{-1}(P)=C_n$. So we can view $C$ as a curve in
$\kc$.
\item[(iii)] For $i=1,\ldots,n$, let \ $\pi_i:\kc_i\to S$ \ be the restriction
of $\pi$. Then \ $\pi_i^{-1}(P)=C$ \ and $\pi_i$ is a flat family of smooth
irreducible projective curves.
\item[(iv)] For every $z\in S\backslash\{P\}$, the components
$\kc_{1,z},\ldots,\kc_{n,z}$ of $\kc_z$ meet transversally and any three
components don't have a common point.
\end{enumerate}
In this case we say also that $\kc$ is a maximal reducible deformation of
\m{C_n}.

Let $L$ be the line bundle on $C$ associated to \m{C_n}. Then for every
\m{z\in S\backslash\{P\}}, any two components \m{\kc_{i,z}}, \m{\kc_{j,z}} of
\m{\kc_z} meet in exactly \m{-\deg(L)} points.
\end{subsub}

\sepprop

\begin{subsub}\label{p2} \rm
Let \m{\kz\subset\kc} be the closure in $\kc$ of the locus of the intersection
points of the components of \m{\pi^{-1}(z)}, \m{z\not=P}. It is a curve in
\m{\kc}.
\end{subsub}

\sepprop

\begin{subsub}\label{p3} \rm
Let \m{I\subset\{1,\ldots,n\}} be a proper subset with $m$ elements, and
\m{\kc_I\subset\kc} the union of the components \m{\kc_i}, \m{i\in I}. Then
the restriction of $\pi$, \m{\pi_I:\kc_I\to S} is a maximal reducible
deformation of \m{C_m}, the inclusion \m{\kc_I\subset\kc} inducing the
inclusion \m{C_m\subset C_n}.
\end{subsub}

\sepprop

\begin{subsub}\label{p4} \rm
Let \ \m{\pi:\kc\to S} \ be a morphism of algebraic varieties satisfying (i),
(iii) and (iv) above, in such a way that the subvarieties \m{C\subset\kc_i}
are identified in $\kc$ and $C$ is the underlying reduced subscheme of
\m{\pi^{-1}(P)}. Then \m{C_n=\pi^{-1}(P)} is a primitive multiple curve of
multiplicity $n$ if and only if for every closed point $x$ of $C$ there exists
an open neighborhood of $x$ in $\kc$ that can be embedded in a smooth variety
of dimension 3 (the proof is the same as that of proposition 4.1.6 of
\cite{dr7}). In this case of course $\pi$ is a maximal reducible deformation
of \m{C_n}.
\end{subsub}

\sepprop

\begin{subsub}\label{p5} Gluings and fragmented deformations -- \rm
For \m{1\leq i\leq n}, let
\m{\pi_i:\kc_i\to S} be a flat family of smooth projective irreducible curves,
with a fixed isomorphism \m{\pi_i^{-1}(P)\simeq C}. A {\em gluing of
\m{\kc_1,\cdots,\kc_n} along $C$} is an algebraic variety $\kd$ such that
\begin{enumerate}
\item[-] for $1\leq i\leq n$, $\kc_i$ is isomorphic to a closed subvariety of
$\kd$, also denoted by $\kc_i$, and $\kd$ is the union of these subvarieties.
\item[-] $\coprod_{1\leq i\leq n}(\kc_i\backslash C)$ is an open subset of
$\kd$.
\item[-] There exists a morphism \ $\pi:\kd\to S$ \ inducing $\pi_i$ on
$\kc_i$, for $1\leq i\leq n$.
\item[-] The subvarieties \ $C=\pi_i^{-1}(P)$ \ of $\kc_i$ coincide in $\kd$.
\end{enumerate}

For example, if \m{\deg(L)=0} in \ref{p1}, i.e. if the irreducible components
of the fibers \m{\pi^{-1}(z)}, \m{z\not=P}, are disjoint, then a maximal
reducible deformation $\kc$ of $C_n$ (which is called a {\em fragmented
deformation} in this case) is a gluing of \m{\kc_1,\cdots,\kc_n} along $C$.

All the gluings of \m{\kc_1,\cdots,\kc_n} along $C$ have the same underlying
Zariski topological space.

Let $\BS{\ka}$ be the {\em initial gluing} of the \m{\kc_i} along $C$. It is an
algebraic variety whose underlying Zarisky topological space is the same as
that of any fragmented deformation with the same components, in particular the
closed points are
\[(\coprod_{i=1}^n\kc_i)/\sim \quad ,\]
where $\sim$ is the equivalence relation: if \m{x\in\kc_i} and
\m{y\in\kc_j}, \m{x\sim y} if and only if \m{x=y}, or if
\m{x\in\kc_{i,P}\simeq C}, \m{y\in\kc_{j,P}\simeq C} and \m{x=y} in $C$.
The structural sheaf is defined by: for every open subset $U$ of $\BS{\ka}$
\[\ko_\BS{\ka}(U) \ = \
\{(\alpha_1,\ldots,\alpha_n)\in\ko_{\kc_1}(U\cap\kc_1)\times
\cdots\times\ko_{\kc_n}(U\cap\kc_n) \ ; \ \alpha_{1\mid 
C}=\cdots=\alpha_{n\mid
C}\} .\]
For every gluing $\kd$ of \m{\kc_1,\cdots,\kc_n}, we have an obvious dominant
morphism \ \m{\BS{\ka}\to\kd}. If follows that the sheaf of rings \m{\ko_\kd}
can be seen as a subsheaf of \m{\ko_{\BS{\ka}}}.

If we consider the maximal reducible deformation $\kc$ of \ref{p1} the
situation is slightly more complicated. We still have a dominant morphism
\ \m{\BS{\ka}\to\kc}. For every \m{x\in C}, \m{\ko_{\kc,x}} can be seen a
subalgebra of \m{\ko_{\BS{\ka},x}}. It contains elements
\m{(\alpha_1,\ldots,\alpha_n)} with the additional property that if \m{1\leq
i<j\leq n} then \m{\alpha_i} and \m{\alpha_j} coincide on \m{\kc_i\cap\kc_j},
which is larger than $C$ if \m{\deg(L)\not=0}. We will see in \ref{frag} that
it is not true that there always exists a fragmented deformation above a
reducible one.
\end{subsub}

\end{sub}

\sepsub

\Ssect{Gluings of families of curves}{arg2b}

Let $S$ be a smooth affine irreducible curve, and \m{P\in S} a closed point.
By a {\em family of curves} parametrized by $S$ we mean a flat projective 
morphism \ \m{\pi:\kc\to S} such that the fibers of $\pi$ are smooth projective 
irreducible curves. The scheme $\kc$ is then a smooth irreducible surface.

\sepprop

\begin{subsub}{\bf Lemma: }\label{lemx1} For every finite subset $\Delta$ of 
$\kc$ there exists an affine open subset \m{U\subset\kc} such that 
\m{\Delta\subset U}. 
\end{subsub}

\begin{proof}
Since $\pi$ is a projective morphism, there exists a projective space \m{\P_N} 
and a closed immersion \ \m{\kc\subset\P_N\times S} \ such that $\pi$ is the 
restriction of the projection \ \m{p_S:\P_N\times S\to S}. Let \ 
\m{p_N:\P_N\times S\to\P_N} \ be the other projection and 
\m{\Delta'=p_N(\Delta)}. There exists an integer \m{k>0} and 
\m{\sigma\in\ko_{P_N}(k)} such that $s$ does not vanish at any point of 
\m{\Delta'}. By \cite{gro1}, prop. 5.5.7 and corr. 1.3.4, the open subset \ 
\m{\kc^\sigma=\{x\in\kc;s(p_N(x))\not=0\}} \ is affine, and it contains 
$\Delta$.
\end{proof}

\sepprop

\begin{subsub}\label{exisfam} Existence of families of curves with prescribed 
sections -- \rm Let $C$ be a smooth projective irreducible curve, and 
\m{x_1,\ldots,x_p} distinct points of $C$.
\end{subsub}

\sepprop

\begin{subsub}\label{lemx2}{\bf Lemma: } For any \m{\eta\in H^1(T_C)}, 
there exists a smooth curve $S$, \m{P\in S}, and a family of curves \ 
\m{\pi:\kc\to S} \ such that
\begin{enumerate}
\item[--] we have \ $\pi^{-1}(P)=C$, and the image of the Koda\" ira-Spencer 
map \ $T_PS\to H^1(T_C)$ \ is $\C\eta$.
\item[--] for \m{1\leq i\leq p}, there exists a section of $\pi$, 
$r_i:S\to\kc$, such that \ $r_i(P)=x_i$.
\end{enumerate}
\end{subsub}

\begin{proof}
If \m{\eta=0} then we can take \ \m{\kc=C\times S} \ with the obvious sections.

Suppose that \m{\eta\not=0}.
There exists a flat family of smooth projective curves \ \m{\theta:\kb\to B} \
parametrized by a smooth variety $B$, such that $\theta$ is a projective 
morphism, and \m{b_0\in B} with
\m{B_{b_0}=\theta^{-1}(b_0)=C}, such that the Koda\" ira-Spencer map \
\m{\omega_{b_0}:TB_{b_0}\to H^1(T_C)} \ is surjective. To obtain $\theta$ one
can use an embedding of $C$ in a projective space \ \m{C\hookrightarrow\P_N} \
such that the degree of $C$ is sufficiently high, take for $B$ a suitable
open subset of the Hilbert scheme of curves in \m{\P_N} with the same genus and
degree as $C$, and for $\kb$ the universal curve (cf. \cite{dr6}, prop. 4.3.1).

Let
\[\kb_p \ = \
\{(\beta_1,\ldots,\beta_p)\in\kb^p \ ; \
\theta(\beta_1)=\cdots=\theta(\beta_p)\} \ . \]
It is a smooth variety and we have
\[T(\kb_p)_{b_0} \ = \ \{(u_1,\ldots,u_p)\in T\kb_{b_0}\times\cdots T\kb_{b_0}
\ ; \ T\theta_{b_0}(u_1)=\cdots=T\theta_{b_0}(u_p)\} \ . \]
Moreover, the projection \m{\kb_d\to B} is a submersion. Let
\[\ov{x} \ = \ (x_1,\ldots,x_p) \ \in \ (\kb_p)_{b_0} \ . \]
Let \m{S\subset\kb_p} be a smooth curve such that \ \m{\ov{x}\in S}.
Let \ \m{f:S\to B} \ be the inclusion \ \m{S\subset\kb_p}, followed by the 
projection \ \m{\kb_p\to B}. Let \ \m{\kc=f^*(\kb)}. It is a family of curves 
parametrized by $S$. It is possible to choose $S$ such that, if \m{u\in 
TS_{\ov{x}}} is a generator, the image in $H^1(T_C)$ (by the Koda\" ira-Spencer 
map of $\kc$) of the projection of $u$ on $TB_{b_0}$ is $\eta$. The family of
curves $\kc$ has all the required properties.
\end{proof}

\sepprop

\begin{subsub}\label{glufam} Gluing of families of curves -- \rm Let $C$ be a 
smooth projective irreducible curve, and \ \m{\pi_1:\kc_1\to S}, 
\m{\pi_2:\kc_2\to S}, families of smooth curves such that \ 
\m{\pi_1^{-1}(P)=\pi_2^{-1}(P)=C}. Let $p$ be a positive integer, and 
\m{x_1,\ldots,x_p} distinct points of $C$. Let \ \m{r_1^i:S\to\kc_1} (resp. 
\m{r_2^i:S\to\kc_2}), \m{1\leq i\leq p}, be sections of \m{\pi_1} (resp. 
\m{\pi_2}), such that \ \m{r_1^i(P)=r_2^i(P)=x_i}. Let
\[\Gamma \ = \ C\cup r_1^1(S)\cup\cdots\cup r_1^p(S) \ \subset \kc_1 \ , \]
which is a closed subvariety of \m{\kc_1}. It is canonically isomorphic to the 
corresponding closed subvariety of \m{\kc_2}:
\[C\cup r_2^1(S)\cup\cdots\cup r_2^p(S) \ . \]
It is the unique scheme obtained by gluing $C$ and $p$ copies of $S$ by 
identifying $x_i$ and $P$ in each copy of $S$, in such a way that all the 
curves $S$ and $C$ are transverse at the intersection points.

According to \cite{ferr}, 4-, 5- th\'eor\`eme 5.4, and lemma \ref{lemx1}, there 
exists a scheme $\kd$, obtained by gluing \m{\kc_1} and \m{\kc_2} along 
$\Gamma$, in such a way that we have a cocartesian diagram
\xmat{\Gamma\flinc[rr]\flinc[dd] & & \kc_1\flinc[dd]\\ \\
\kc_2\flinc[rr] & & \kd=\kc_1\sqcup_\Gamma\kc_2}
If \m{x\in\Gamma}, \m{\ko_{\kd,x}} is the ring of pairs \ \m{(\phi_1,\phi_2)\in
\ko_{\kc_1,x}\times\ko_{\kc_2,x}} \ such that \ 
\m{\phi_{1|\Gamma}=\phi_{2|\Gamma}}. Moreover the morphisms \m{\pi_1}, \m{\pi_2}
can be extended to a morphism \ \m{\pi:\kd\to S}.

Let \m{\ko_{\kc_i}(1)}, \m{i=1,2}, be ample line bundles on \m{\kc_i}, such that
\ \m{\ko_{\kc_1}(1)_{|C}\simeq\ko_{\kc_2}(1)_{|C}}, and that 
\m{\ko_{\kc_i}(1)_{r_i^j(S)}} is trivial for \m{1\leq j\leq p}. Then we have \ 
\m{\ko_{\kc_1}(1)_{|\Gamma}\simeq\ko_{\kc_2}(1)_{|\Gamma}}. By the same 
argument as in \cite{ferr}, 6.3, we can define a line bundle \m{\ko_\kd(1)} on 
$\kd$ extending \m{\ko_{\kc_1}(1)} and \m{\ko_{\kc_2}(1)}. By \cite{gro2}, 
2.6.2, applied to \ \m{\pi:\kd\to S} \ and \ 
\m{(\pi_1,\pi_2):\kc_1\sqcup\kc_2\to S}, \m{\ko_\kd(1)} is ample. Hence $\kd$ 
is a quasiprojective variety. According to \cite{dr7} (cf. 2.3), \m{\pi^{-1}(P)} 
is a primitive double curve with associated smooth curve $S$, and $\kd$ is a 
maximal reducible deformation of \m{\pi^{-1}(P)}.
\end{subsub}

\end{sub}

\sepsub

\Ssect{Smoothing of reducible deformations of primitive multiple curves}{arg2}

Let $X$ be an algebraic variety. We say that $X$ is {\em smoothable} if there
exists a flat morphism \ \m{\psi:\kx\to T}, where $\kx$ and $T$ are algebraic
varieties, such that $T$ is integral, there exists \m{t_0\in T} such that
\m{\psi^{-1}(t_0)\simeq X}, and if \m{s\not=t\in T}, then \m{\psi^{-1}(s)} is
smooth.

\sepprop

\begin{subsub}{\bf Proposition : }\label{pro_def} Let $D$ be a projective
primitive multiple curve. If there exists a maximal reducible deformation of
$D$, then $D$ is smoothable.
\end{subsub}

\begin{proof}
Let \ \m{\pi:\kc\to S} \ be a maximal reducible deformation of $D$, and
\m{s_0\in S} such that \m{D=\kc_{s_0}}. Let \ \m{i:\kc\to\P_N} \ be an
embedding of $\kc$ in a projective space. We may assume that for every \m{s\in
S} we have \ \m{h^1(\ko_{\kc_s}(1))=0}, and the exact sequence
\[0\lra\ko_{\P_N}\lra\ko_{\P_N}(1)\ot H^0(\ko_{\P_N}(1))^*\lra T\P_N\lra 0\]
implies that we have also \ \m{h^1(T\P_{N\mid\kc_s})=0}. Let $\kz$ be
the component of the Hilbert scheme of curves in \m{\P_N} containing $D$.
Then by the local structure of \m{D} there is a canonical surjective morphism
\xmat{T\P_N\flon[r] & \HHom(\ki_D,\ko_D)}
(where \m{\ki_D} denotes the ideal sheaf of $D$ in \m{\P_N}). It follows that \
\m{H^1(\HHom(\ki_D,\ko_D))=\nsp}. The same holds if we replace $D$ with
\m{\kc_s}, for every \m{s\in S}. It follows that $\kz$ is smooth at all the
fibers of $\pi$. Let \m{\kz^0\subset\kz} be the open subset of smooth points.
We have thus \ \m{\kc_s\in\kz^0} for every \m{s\in S}.

It is clear that all the singularities of the fibers \m{\kc_s}, \m{s\in
S\backslash\{s_0\}}, are smoothable. It follows from proposition 29.9 of
\cite{ha2} that the fibers \m{\kc_s}, \m{s\in S\backslash\{s_0\}}, are
smoothable in \m{\P_N}. In particular some points in \m{\kz^0} are smooth
curves. It follows that \ \m{D=\kc_{s_0}} \ is smoothable.
\end{proof}

\end{sub}

\sepsec

\section{Maximal reducible deformations of primitive double
curves}\label{Maxred2}

In this chapter $S$ denotes a smooth curve, and \m{P\in S}. Let
\m{t\in\ko_{S,P}} be a generator of the maximal ideal of $P$. We can suppose
that $t$ is defined on the whole of $S$, and that the ideal sheaf of $P$ in
$S$ is generated by $t$.

\sepsub

\Ssect{Properties of maximal reducible deformations of primitive double
curves}{M2prop}

Let \m{C_2} be a primitive double curve, with underlying projective smooth
curve $C$ and associated line bundle $L$ on $C$.  We suppose that \ 
\m{\deg(L)<0}. Let \ \m{\pi:\kc\to S} \ be a maximal reducible deformation of 
\m{C_2}, and \m{P\in S} such that \ \m{\pi^{-1}(P)=C_2}. Then $\kc$ has two 
irreducible components \m{\kc_1}, \m{\kc_2} which are flat families of smooth 
irreducible curves parametrized by $S$. If \m{z\in S\backslash\{P\}}, the two 
components \m{\kc_{1,z}} and \m{\kc_{2,z}} of \m{\kc_z} intersect transversally 
in \m{-\deg(L)} points.

Let \m{\kz\subset\kc} be the closure in $\kc$ of the locus of the intersection
points of the components of \m{\pi^{-1}(z)}, \m{z\not=P}. Since $S$ is a curve,
$\kz$ is a curve of \m{\kc_1} and \m{\kc_2}. It intersects $C$ in a finite
number of points. If \m{x\in C}, let \m{r_x} be the number of
branches of $\kz$ at $x$ and \m{s_x} the sum of the multiplicities of the
intersections of these branches with $C$, so that we have \ \m{r_x\leq s_x},
with equality (for all such $x$) if and only if all the branches intersect 
transversally with $C$. Moreover, since for every \ \m{z\in S\backslash\{P\}},
\m{\kc_{1z}\cap\kc_{2z}} consists of \m{-\deg(L)} distinct points, we have \
\m{\sigg_{x\in\kz\cap C}r_x=-\deg(L)}.

For \m{i=1,2}, let \m{\pi_i:\kc_i\to S} be
the restriction of $\pi$. We will also denote \m{\pi^*t} by $\pi$, and
\m{\pi_i^*t} by \m{\pi_i}. So we have
\m{\pi=(\pi_1,\pi_2)\in\ko_\kc(-C)}.

\sepprop

\begin{subsub}\label{th1}{\bf Theorem:  1 -- } Let \m{x\in C}. Then there
exists an unique integer \m{p>0} such that \m{\ki_{C,x}/\span{(\pi_1,\pi_2)}}
is generated by the image of \m{(\pi_1^p\lambda,0)}, for some
\m{\lambda\in\ko_{\kc_1,x}} not divisible by \m{\pi_1}. This integer does not
depend on $x$.

{\bf 2 -- } $\lambda$ is unique up to multiplication by an invertible element
of \m{\ko_{\kc_1,x}}, and \m{(\pi_1^p\lambda,0)} is a generator of the ideal
\m{\ki_{\kc_1,\kc,x}} of \m{\kc_1} in $\kc$.

{\bf 3 -- } There are only a finite number of points \m{x\in C} such that
$\lambda$ is not invertible.

{\bf 4 -- } Let \m{m_x} be the multiplicity of \m{\lambda_{\mid
C}\in\ko_{C,x}}. Then we have \ \m{m_x>0} \ if and only if \ \m{x\in\kz\cap
C}, and in this case we have \ \m{m_x=r_x=s_x}, and the branches of $\kz$ at
$x$
intersect transversally with $C$. Moreover
\[L \ \simeq \ \ko_C(-\sigg_{x\in \kz\cap C}r_xx) \ \simeq \ \ki_{\kz\cap C,C}
\ . \]
\end{subsub}

\begin{proof} the proof of 1- is similar to the proof of proposition
4.2.1, 1- of \cite{dr7}. Let \m{x\in C} and \m{u=(\pi_1\alpha,\pi_2\beta)} whose
image is a generator of \m{\ki_C/\span{(\pi_1,\pi_2)}} at $x$
(\m{\ki_C/\span{(\pi_1,\pi_2)}} is a locally free sheaf of rank 1 of
\m{\ko_C}-modules). Let \m{\beta_0\in\ko_{\kc_1,x}} be such that
\m{(\beta_0,\beta)\in\ko_{\kc,x}}. Then the image of
\[u-(\pi_1,\pi_2)(\beta_0,\beta)=(\pi_1(\alpha-\beta_0),0)\]
is also a generator of \m{\ki_C/\span{(\pi_1,\pi_2)}} at $x$. We can write it
\m{(\pi_1^p\lambda,0)}, where $\lambda$ is not a multiple of $\pi_1$.

As in proposition 4.2.1 of \cite{dr7} it is easy to see that if
\m{(\pi_1^q\mu,0)\in\ko_{\kc,x}}, with $\mu$ not divisible by $\pi_1$ and
\m{q>0}, and if $k$ is a positive integer, then we can write
\[(\pi_1^q\mu,0) \ = \ \gamma.(\pi_1^p\lambda,0)+\delta.(\pi_1^k,\pi_2^k) \ ,
\]
with \m{\gamma,\delta\in\ko_{\kc,x}}. From this it follows that \m{q\geq p},
so $p$ is unique. Since the image of \m{(\pi_1^p\lambda,0)} is also a
generator of \m{\ki_C/\span{(\pi_1,\pi_2)}} in a neigborhood of $x$, $p$ does
not depend on $x$ and 1- is proved, as well as 3-.

Suppose that the image of \m{(\pi_1^p\lambda',0)} is also a generator
of \m{\ki_C/\span{(\pi_1,\pi_2)}} at $x$. Then for every positive integer $k$
we can write
\[(\pi_1^p\lambda',0) \ = \ (\gamma_1,\gamma_2)(\pi_1^p\lambda,0)+
(\delta_1,0).(\pi_1^{p+k},\pi_2^{p+k}) \ , \]
with \m{(\gamma_1,\gamma_2),(\delta_1,0)\in\ko_{\kc,x}}, whence \
\m{\lambda'=\gamma_1\lambda+\pi_1^k\delta_1}. We have then
\[\lambda' \ \in \ \bigcap_{k>0}((\lambda)+(\pi_1)^{p+k}\]
in \m{\ko_{\kc_1,x}}, and the latter is equal to \m{(\lambda)} according to
\cite{sa-za}, vol. II, chap. VIII, theorem 9. So we can write \
\m{\lambda'=\beta\lambda} \ with \m{\beta\in\ko_{\kc_1,x}}. In the same way,
\m{\lambda=\beta'\lambda'} \ with \m{\beta'\in\ko_{\kc_1,x}}. So
\m{\beta\beta'=1} and $\beta$ is invertible, which proves the first assertion
2-. The proof of the second assertion is similar.

Suppose that \m{x\in\kz\cap C}. Since \ \m{(\pi_1^p\lambda,0)\in\ko_{\kc,x}},
we have \ \m{\lambda_{\mid\kz\cap C}=0}, hence \ \m{m_x\geq s_x} .

It follows from 1-, 2-, 3- that \m{\ki_C/\span{(\pi_1,\pi_2)}} (which is
isomorphic to $L$) can be viewed as a subsheaf of \
\m{(\pi_1^p)/(\pi_1^{p+1})\simeq\ko_C}, and this subsheaf is exactly \
\m{\ko_C(-\sigg_{x\in \kz\cap C}m_xx)}. Hence \ \m{-\deg(L)=\sigg_{x\in C}m_x}.
So
\[-\deg(L)=\sigg_{x\in C}m_x \ \geq \sigg_{x\in\kz\cap C}s_x \ \geq \
\sigg_{x\in\kz\cap C}r_x=-\deg(L) \ , \]
hence all the preceding inequalities are equalities and 4- follows.
\end{proof}

\sepprop

Similarly there exists an unique integer \m{q>0} such that for every \m{x\in
C}, \m{\ki_{C,x}/\span{(\pi_1,\pi_2)}} is generated by the image of
\m{(0,\pi_2^q\mu)}, for some \m{\mu\in\ko_{\kc_2,x}} not divisible by
\m{\pi_2} and unique up to multiplication by an invertible element of
\m{\ko_{\kc_2,x}}.

\sepprop



\begin{subsub}\label{coro1}{\bf Corollary: } Let \m{x\in C} and \
\m{(\pi_1^p\lambda,0)\in\ko_{\kc,x}} (resp. \m{(0,\pi_2^p\mu)\in\ko_{\kc,x}})
be such that its image is a generator of \ \m{\ki_{C,x}/\span{(\pi_1,\pi_2)}}.
Let \ \m{\theta\in\ko_{\kc_2,x}} \ such that \
\m{(\lambda,\theta)\in\ko_{\kc,x}}. Then

{\bf 1 -- } There exists \m{\gamma\in\ko_{\kc_2,x}} invertible such that
\m{\theta=\gamma\mu}.

{\bf 2 -- } The ideal of $\kz$ in \m{\kc_{1,x}} (resp. \m{\kc_{2,x}}) id
generated by $\lambda$ (resp. $\mu$).
\end{subsub}

\begin{proof} Since the image of \
\m{\pi^p(\lambda,\theta)-(\pi_1^p\lambda,0)=(0,-\pi_2^p\theta)} \ is also a
generator of \ \m{\ki_{C,x}/\span{(\pi_1,\pi_2)}}, 1- follows from
theorem \ref{th1}, 2- (for the other coordinate). The second assertion follows
from theorem \ref{th1}, 4-.
\end{proof}

\sepprop

\begin{subsub}\label{prop1}{\bf Proposition: 1 -- } There exists a canonical
isomorphism
\[\Phi:\ko_{\kc_1}/(\pi_1^p\ki_{\kz,\kc_1})\lra
\ko_{\kc_2}/(\pi_2^q\ki_{\kz,\kc_2}) \]
such that for every \m{x\in C} and \ \m{(\alpha_1,\alpha_2)\in\ko_{\kc_1,x}
\times\ko_{\kc_2,x}}, if \m{[\alpha_1]}, \m{[\alpha_2]} denote the images of
\m{\alpha_1}, \m{\alpha_2} in \
\m{\ko_{\kc_1,x}/(\pi_1^p\ki_{\kz,\kc_1,x})},
\m{\ko_{\kc_2,x}/(\pi_1^p\ki_{\kz,\kc_2,x})} \ respectively, we have
\ \m{(\alpha_1,\alpha_2)\in\ko_{\kc,x}} \ if and only if \ \m{\Phi([\alpha_1])
=[\alpha_2]}. For every \ \m{\alpha\in\ko_{\kc,x}}, we have \
\m{\Phi([\alpha])_{\mid C}=[\alpha]_{\mid C}}, and \ \m{\Phi([\pi_1])=[\pi_2]}.

{\bf 2 -- } We have \m{q=p}.
\end{subsub}

\begin{proof}
The first assertion is an easy consequence of the second statement of the
preceding theorem. To prove 2-, let \m{x\in C\backslash\kz\cap C}. Then
\m{(\pi_1^p\lambda,0)\in\ko_{\kc,x}} for some \m{\lambda\in\ko_{\kc_1,x}} such
that \m{\lambda_{\mid C}\not=0}. There exists \m{\mu\in\ko_{\kc_2,x}} such
that \m{(\lambda,\mu)\in\ko_{\kc,x}}, and \m{\mu_{\mid C}=\lambda_{\mid
C}\not=0}. We have \
\m{(0,\mu\pi_2^p)=\pi^p(\lambda,\mu)-(\pi_1^p\lambda,0)\in\ko_{\kc,x}}, so
\m{q\leq p} \ by theorem \ref{th1}, 2-. Similarly \m{p\leq q}, hence \m{p=q}.
\end{proof}

\sepprop

\begin{subsub}\label{prop1b}{\bf Proposition: } Let \m{J_p\subset\ko_\kc}
be the ideal sheaf consisting of pairs \m{(u_1,u_2)} such that \m{u_1} is a
multiple of \m{\pi_1^{p+1}} and \m{u_2} a multiple of \m{\pi_2^{p+1}}. Then we
have \ \m{J_p\subset(\pi)}.
\end{subsub}

\begin{proof}
Let \m{a_1\in\ko_{\kc_1,x}, a_2\in\kc_{2,x}} be such that \
\m{\Phi([a_1\pi_1^{p+1}])=[a_2\pi_2^{p+1}]}. We must prove \Nligne that \
\m{\Phi([a_1\pi_1^p])=[a_2\pi_2^p]}. Let \m{b\in\ko_{\kc_2,x}} be such that \
\m{\Phi([a_1])=[b]}. Then we have \Nligne
\m{\Phi([a_1\pi_1^{p+1}])=[b\pi_2^{p+1}]},
hence \m{a_2\pi_2^{p+1}-b\pi_2^{p+1}} is a multiple of \m{\mu\pi_2^p} :
\[a_2\pi_2^{p+1}-b\pi_2^{p+1} \ = \beta\mu\pi_2^p\]
for some \m{\beta\in\ko_{\kc_{2,x}}}. It follows that $\beta$ is a multiple of
\m{\pi_2} : \m{\beta=\alpha\pi_2}, and \ \m{a_2-b=\alpha\mu}. Hence
\begin{eqnarray*}
\Phi([a_1\pi_1^p]) & = & [b\pi_2^p]\\
& = & [a_2\pi_2^p-\alpha\mu\pi_2^p]\\
& = & [a_2\pi_2^p] \ .
\end{eqnarray*}
\end{proof}

\sepprop

It follows that for the associated primitive multiple curve \m{C_2} we have
\[\ko_{C_2} \ = \ \ko_\kc/(\pi) \ = \ \big(\ko_\kc/J_p\big)/(\pi) \ . \]

\sepprop

\begin{subsub}\label{local} Localization -- \rm Let
\[\ov{\Phi}: \ko_{\kc_1}/\big((\pi_1^p\ki_{\kz,\kc_1})+(\pi_1^{p+1})\big)\lra
\ko_{\kc_2}/\big(\pi_2^q\ki_{\kz,\kc_2})+(\pi_1^{p+1})\big)\]
be the isomorphism induced by $\Phi$. Let \m{x\in C}. For \m{i=1} or 2, and
\m{\alpha_i\in\ko_{\kc_i,x}/(\pi_i^{p+1})}, let \m{\epsilon(\alpha_i)} denote
the image of \m{\alpha_i} in
\m{\big[\ko_{\kc_i}/\big((\pi_i^p\ki_{\kz,\kc_i})+(\pi_1^{p+1})\big)\big]_x}.
Let
\[\ka_x \ = \ \{(\alpha_1,\alpha_2)\in\ko_{\kc_1,x}/(\pi_1^{p+1})\times
\ko_{\kc_2,x}/(\pi_2^{p+1});\ov{\Phi}(\epsilon(\alpha_1))=\epsilon(\alpha_2)
\} \ . \]
\end{subsub}

\sepprop

\begin{subsub}\label{prop1c}{\bf Proposition: } The natural projection \
\m{\ko_{\kc,x}\to\ka_x} \ induces an isomorphism
\[\ko_{\kc,x}/(\pi) \ \simeq \ \ka_x/(\pi) \ . \]
\end{subsub}

\begin{proof}
Let \ \m{\theta:\ko_{\kc,x}/(\pi)\to\ka_x/(\pi)} \ be the natural morphism.

First we prove that $\theta$ is surjective. Let \m{(\alpha_1,\alpha_2)\in\ka_x}.
Let \m{a_i\in\ko_{\kc_i,x}} be over \m{\alpha_i}. Then, since \
\m{\ov{\Phi}(\epsilon(a_1))=\epsilon(a_2)}, we can write
\[\Phi([a_1]) \ = \ [a_2]+[\beta_2]\pi_2^{p+1} \ , \]
for some \m{\beta_2\in\ko_{\kc_2,x}}. Let \m{a'_2=a_2+\pi_2^{p+1}\beta_2}. Then
the image of \m{(a_1,a'_2)} in \m{\ka_x} is \m{(\alpha_1,\alpha_2)}, and
clearly \m{(a_1,a'_2)\in\ko_{\kc_x}} and its image in \m{\ka_x} is
\m{(\alpha_1,\alpha_2)}.

Now we prove that $\theta$ is injective. Let \m{u=(u_1,u_2)\in\ko_{\kc,x}},
$\ov{u}$ its image in \m{\ko_{\kc,x}/(\pi)}, and suppose that
\m{\theta(\ov{u})=0}. Let $v$ be the image of $u$ in \m{\ka_x}. Then we can
write \m{v=\pi w}, for some \m{w\in\ka_x}. Since $\theta$ is surjective, we can
find \m{w'\in\ko_{\kc,x}} over $w$. Then \m{u'=u-\pi w'} is of the form
\m{(\pi_1^{p+1}\alpha_1,\pi_2^{p+1}\alpha_2)}, for some
\m{\alpha_i\in\ko_{\kc_i,x}}. From proposition \ref{prop1b}, \m{u'} is a
multiple of $\pi$. Hence $u$ is a multiple of $\pi$ and \m{\ov{u}=0}.
\end{proof}

\end{sub}

\sepsub

\Ssect{Construction of maximal reducible deformations of primitive double
curves}{const2}

We keep the notations of \ref{M2prop} and we suppose that \m{p=1} (we don't 
need to study the case \m{p>1} to prove theorem \ref{th2}, (ii)).

\sepprop

\begin{subsub}\label{justif} Maximal reducible deformations of primitive double
curves and gluings of families of curves -- \rm We consider the maximal 
reducible deformation \m{\pi:\kc\to S} of \ref{M2prop} (with \m{p=1}).
It follows from the description of the local rings of the closed points of 
$\kc$ 
given in \ref{M2prop} that $\kc$ is a {\em gluing} of $\kc_1$ and \m{\kc_2} 
along the closed subvariety \m{\Gamma=\kz\cup C}. This means that we have a 
cocartesian diagram
\xmat{\Gamma\flinc[rr]\flinc[dd] & & \kc_1\flinc[dd]\\ \\
\kc_2\flinc[rr] & & \kc=\kc_1\sqcup_\Gamma\kc_2}
and that for every closed point \m{x\in\Gamma}, \m{\ko_{\kc,x}} is the ring of 
pairs \ \m{(\phi_1,\phi_2)\in \ko_{\kc_1,x}\times\ko_{\kc_2,x}} \ such that \ 
\m{\phi_{1|\Gamma}=\phi_{2|\Gamma}}. This comes also from the fact that \ 
\m{\Gamma=\kc_1\cap\kc_2} as schemes (of course this is not true if \m{p>1}).
\end{subsub}

\sepprop

\begin{subsub}\label{const2_1} The cocycles defining \m{C_2} -- \rm
For \m{i=1,2}, let \m{\kc_i^{(2)}} denote the infinitesimal neighborhood of
order 2 of $C$ in \m{\kc_i} (i.e. \
\m{\ko_{\kc^{(2)}_i}=\ko_{\kc_i}/(\pi_i^2)}). It is a primitive double curve.
Let \m{(U_j)_{j\in J}} be a finite open affine cover of $C$. Then for every
\m{j\in J}, the restriction \m{U_j^{(2)}} of \m{\kc_i^{(2)}} to \m{U_j} is
trivial, i.e. there is an isomorphism \ \m{U_j^{(2)}\simeq U_j\times Z_2}
(where \ \m{Z_2=\spec(\C[t])/(t^2)}) inducing the identity on \m{U_j}. We can
suppose that \m{U_j} contains at most one point of \m{\kz\cap C}, and that each
point in \m{\kz\cap C} is contained in only one \m{U_j}. We suppose
also that for every distinct \m{j,k\in J}, \m{\omega_{C|U_{jk}}} is trivial,
generated by \m{dx_{jk}}. Then \m{\kc_i^{(2)}} can be constructed by a cocycle
\m{(\mu_{jk}^{(i)})_{j,k\in J}}, \m{\mu_{jk}^{(i)}\in\ko_C(U_{jk})}, as in
\ref{dblcst}: \m{\kc_i^{(2)}} is obtained by gluing the \ \m{U_j^{(2)}\simeq
U_j\times Z_2} \ with the automorphisms of \ \m{U_{jk}\times Z_2} \ defined by
the matrices \ \m{\dsp\begin{pmatrix}1 & 0 \\
\mu_{jk}^{(i)}\frac{\partial}{\partial x_{jk}} & 1\end{pmatrix}} \ .

If \m{U_j} does not contain any point of \m{\kz\cap C}, let \m{r_j^ {(i)}=1}. If
\m{U_j} contains a point of \m{\kz\cap C}, this point $x$ is unique, and in
this case let \m{r_j^{(i)}\in\ko_C(U_j)[t]/(t^2)} be an equation
of \m{\kz\cap C} in the open subset \m{U_j\times Z_2} of \m{\kc^{(i)}_2}. If
\m{\tau\in\ko_C(U_j)} is an equation of $x$, we can write (by theorem \ref{th1})
\[r_j^{(i)} \ = \ \tau^{r_x}+t.h_j^{(i)} \ , \]
where \m{h_j^{(i)}} vanishes to order \m{\geq r_x-1} at $x$. It follows that we
can suppose that \ \m{r_j^{(1)}=r_j^{(2)}}. More precisely, there exists an
automorphism $\sigma$ of \m{\ko_C(U_j)[t]/(t^2)} such that the induced
automorphism of \m{\ko_C(U_j)} is the identity, and sending \m{r_j^{(1)}} to
\m{r_j^{(2)}} (this comes from the description of these automorphisms in
\ref{constpmc}). Let \ \m{r_j=r_j^{(1)}=r_j^{(2)}}. Note that by theorem
\ref{th1} we can assume that \m{r_j} is a product
\[r_j \ = \ \prod_{m=1}^{r_x}(\tau+\lambda_mt) \ , \]
where \m{\lambda_m\in\C} for \m{1\leq m\leq r_x}, and the \m{\lambda_m} are
distinct. Let \ \m{\rho_j=r_{j|C}}.

Let \m{\kc^{(2)}} be the scheme corresponding to the sheaf of algebras $\ka$
(cf. \ref{local}). We have \ \m{(\kc^{(2)})_{red}=C}. For every \m{j\in J}, we
can view \m{\ka(U_j)} as the algebra of pairs \m{(a+bt,a+(b+\rho_j\beta)t)},
with \m{a,b,\beta\in\ko_C(U_j)} and the rule \m{t^2=0}. These sets must then be
glued to build $\ka$, using the automorphisms, for distinct \m{j,k\in J}
\begin{equation}\label{equ001}\xymatrix@R=5pt{\ka(U_{jk})\ar[r] & \ka(U_{jk})\\
(a+bt,a+(b+\rho_j\beta)t)\fmaps[r] & \big(a+(\mu_{jk}^{(1)}\frac{\partial
a}{\partial x_{jk}}+b)t,a+(\mu_{jk}^{(2)}\frac{\partial a}{\partial x_{jk}}+b
+\rho_j\beta)t\big)} \ . \end{equation}

Now we can also describe the double primitive curve \m{C_2} with a cocycle,
unsing the fact that \ \m{\ko_{C_2}=\ko_\kc/(\pi)=\ka/(\pi)} (cf. prop.
\ref{prop1c}). If \m{u\in\ka(U_j)}, let \m{[u]} denote its image in
\m{\ka(U_j)/(\pi)}. In \m{\ka(U_j)} we have \ \m{\pi=(t,t)}, hence we have an
isomorphism
\[\xymatrix@R=5pt{(\ka/(\pi))(U_j)\ar[r] & \ko_C(U_j)[z]/(z^2)\\
(a+bt,a+(b+\rho_j\beta)t)\fmaps[r] & a+\beta z \ . }\]

Hence the automorphism of \m{\ko_C(U_{jk})[z]/(z^2)} defining \m{C_2} is given
by
\[a+bz\longmapsto a +\big(\frac{\mu_{jk}^{(2)}-\mu_{jk}^{(1)}}{\rho_j}
\frac{\partial a}{\partial x_{jk}}+\frac{\rho_j}{\rho_k}b\big)z \]
(note that \m{\rho_j} and \m{\rho_k} are invertible on \m{U_{jk}}). It follows 
that \m{C_2} is defined by the cocycles 
\m{\dsp\left(\frac{\rho_j}{\rho_k}\right)} and
\m{\dsp\left(\frac{\mu_{jk}^{(2)}-\mu_{jk}^{(1)}}{\rho_j}\right)}. It is easy
to verify that the first one defines \ \m{L=\ko_C(-\sigg_{x\in\kz\cap C}r_xx)}.
\end{subsub}

\sepprop

\begin{subsub}\label{loc_red_def} Local reducible deformations -- \rm Let
\m{\kd_1}, \m{\kd_2} be primitive double curves with associated smooth curve
$C$ and associated line bundle \m{\ko_C}. Suppose that \m{\kd_i}, \m{i=1,2}, is
defined by the cocycle \m{(\rho^{(i)}_{jk}\frac{\partial}{\partial x_{jk}})}
(it is not necessary to use another cocycle since \m{L=\ko_C}). For every
\m{j\in J} let \m{r_j\in\ko_C(U_j)[t]/(t^2)} such that
\begin{enumerate}
\item[--] $r_j$ vanishes in at most one point in $U_j$, and not on $U_j\cap
U_k$ if $k\in J\backslash\{j\}$.
\item[--] if $r_j$ vanishes at $x_j\in U_j$ and $\tau$ is a generator of the
ideal of $x$ in $\ko_C(U_j)$, then \ $r_j=\prod_{m=1}^{p_j}(\tau+\lambda_mt)$,
for some integer $p_j\geq 1$, where the $\lambda_m$ are distinct scalars.
\end{enumerate}
Let \m{J'\subset J} be the set of points $j$ such that \m{r_j} vanishes at some
point of \m{U_j}, and $Z$ be the divisor \m{\sigg_{j\in J'}p_jx_j}.

Then we can define from these data, using a sheaf of algebras $\ka$ as in
\ref{const2_1}, a scheme $\kd$ such that \m{\kd_{red}=C}, and
a global section $\pi$ of \m{\ko_\kd}. We simply define \m{\ka(U_j)} as the
algebra of pairs \m{(a+bt,a+(b+r_j\beta)t)}, with \m{a,b,\beta\in\ko_C(U_j)}
and the rule \m{t^2=0}. We then glue the schemes \m{\spec(\ka(U_j))} using
automorphisms similar to (\ref{equ001}). Now the sheaf \m{\ka/(\pi)} if the
structure sheaf of a double primitive curve \m{C_2} with associated smooth
curve $C$ and associated line bundle \ \m{L=\ko_C(-Z)}.

We call the scheme $\kd$ a {\em local reducible deformation} of the primitive
double curve \m{C_2}. The reducible deformations $\kc$ defined previously can
also be called {\em global} reducible deformations. It is clear from 3.2.2 that 
a global reducible deformation of \m{C_2} induces a local one.

\sepprop

\begin{subsub}\label{rem1}{\bf Remark: } \rm It is possible to define local 
reducible deformations in the case \m{p>1}. But in this case we would need to 
use schemes of the form \ \m{U_i\times Z_{p+1}} instead of \m{U_i\times Z_2}.
\end{subsub}

\sepprop

\begin{subsub}\label{prop1d}{\bf Proposition: }
Let $D$ be a primitive double curve with associated smooth curve $C$ and
associated line bundle \m{L'} on $C$. Suppose that \ \m{h^0({L'}^*)\not=0}.
Then there exists a local reducible deformation of $D$.
\end{subsub}

\begin{proof} The hypothesis \ \m{h^0({L'}^*)\not=0} \ means that \m{L'} is an
ideal sheaf, so we can write \ \m{L'=\ko_C(-\sigg_{i=1}^pn_ix_i)}, where
\m{x_1,\ldots,x_p} are distinct points of $C$ and \m{n_1,\ldots,n_p} positive
integers. We can choose an open affine cover \m{(U_i)_{i\in I}} of $C$ such
that any \m{U_i} contains at most one of the points \m{x_j}, and that every
point \m{x_j} is contained in only one \m{U_i}. We take now
\m{r_j\in\ko_C(U_j)[t]/(t^2)} such that
\begin{enumerate}
\item[--] $r_j$ is invertible if $U_j$ contains no point $x_i$.
\item[--] If $x_i\in U_j$, and $\tau\in\ko_C(U_j)$ a generator of the ideal of
$x_i$, then $r_j$ is of the form \ $r_j=\prod_{1\leq m\leq
n_i}(\tau+\lambda_mt)$, where $\lambda_1,\ldots,\lambda_m\in\C^*$.
\end{enumerate}
For \m{1\leq i\leq p}, let \m{\rho=r_{i|C}}. Then
\m{(\dsp\frac{\rho_j}{\rho_k})} is a cocycle of invertible functions (with
respect to \m{(U_{jk})}) which defines \m{L'}. Suppose that $D$ comes from \
\m{\sigma\in H^1(T_C\ot L')}, defined by the cocycle \m{(\nu_{jk})}. Now we
define two primitive double curves with associated smooth curve $C$ and
associated line bundle \m{\ko_C}:
\begin{enumerate}
\item[--] The first one in the trivial primitive double curve, defined by the
cocycle $(0)$.
\item[--] The second one is defined by the cocycle $(\rho_j\nu_{jk})$.
\end{enumerate}
It is clear from \ref{const2_1} that these data define a local reducible
deformation of $D$.
\end{proof}

\sepprop

We have seen that a reducible deformation of \m{C_2} induces a local reducible
deformation of it. The converse is true (for suitable $L$):
\end{subsub}

\sepprop

\begin{subsub}\label{prop1e}{\bf Proposition: } Let $D$ be a primitive double 
curve, with associated smooth curve $C$ and associated line bundle $L$. 
Suppose that there exist distinct points \m{P_1,\ldots,P_m}of $C$, such that
\[L \ \simeq \ \ko_C(-P_1-\cdots-P_m) \ . \]
Then there exists a global maximal reducible deformation of $D$.
\end{subsub}

\begin{proof}
We begin as in \ref{prop1d} with a local reducible deformation of $D$. We 
have to show that it can be extended to a global reducible deformation of $D$. 
This is an immediate consequence of lemma 2.4.3 (for the construction of 
adequate \m{\kc_1} and \m{\kc_2}) and 2.4.4 (to glue them and make $\kc$).
\end{proof}

\sepprop

The preceeding results can then be summarized in

\sepprop

\begin{subsub}\label{th2}{\bf Theorem: } Let $D$ be a primitive double curve
with associated smooth curve $C$ and associated line bundle $L$ on $C$. Then

{\em (i)} There exists a local maximal reducible deformation of $D$ if and only 
if \ \m{h^0(L^*)\not=0}.

{\em (ii)} If there exists a divisor \m{\Delta=P_1+\cdots+P_m} of $C$, with 
distinct points \m{P_1,\ldots,P_m} of $C$, such that \ 
\m{L\simeq\ko_C(-\Delta)}. Then there exists a global maximal reducible 
deformation of $D$.
\end{subsub}

\sepprop

\begin{subsub}\label{exa}{\bf Example : } \rm Let $\Delta$ be a divisor as in 
theorem \ref{th2}, and \ \m{L=\ko_C(-\Delta)}. Let \ \m{\kc_1=\kc_2=C\times\C} 
and
\[\Gamma \ = \ C\cup(\{P_1\}\times\C)\cup\cdots\cup(\{P_m\}\times\C) \ . \]
This curve is a closed subvariety of \m{\kc_1} and \m{\kc_2}. We can glue 
\m{\kc_1} and \m{\kc_2} along $\Gamma$ (using 2.4.4), and we obtain a maximal 
reducible deformation $\kc$ of the trivial double curve associated to $C$ and 
$L$.

Another way to obtain $\kc$ is as follows: let \m{s_0} be the zero section of 
\m{L^*}, and\ \m{s_1\in H^0(L^*)} \ vanishing at \m{P_1,\ldots,P_m}. We can 
view $C$ as embedded in the surface \m{L^*} via \m{s_0}. Then we can realize 
\m{\kc_1} and \m{\kc_2} as closed subvarieties of \ \m{L^*\times\C} \ in the 
following way : \m{\kc_1=C\times\C}, and \m{\kc_2} is over \m{t\in\C} the image 
of the section \m{s_0+ts_1}. Then we have \ \m{\kc=\kc_1\cup\kc_2}.
\end{subsub}

\end{sub}

\sepsub

\Ssect{Smoothing of primitive double curves}{smo_pri_dou}

The following result follows immediately from proposition \ref{pro_def} and
theorem \ref{th2}: let $D$ be a primitive double curve
with associated smooth curve $C$ and associated line bundle $L$ on $C$.
If there exists a divisor \m{P_1+\cdots+P_m} of $C$, with distinct points 
\m{P_1,\ldots,P_m} of $C$, such that \ \m{L\simeq\ko_C(-P_1-\cdots-P_m)}, then 
$D$ is smoothable.
Of course this is also a consequence of \cite{gon} (cf. the Introduction), 
since we have also \ \m{h^0(L^{-2})\not=0}.

\end{sub}

\sepsec

\section{Maximal reducible deformations in the general case}\label{Maxredn}

\Ssect{Properties of maximal reducible deformations of primitive multiple
curves}{Mnprop}

Let $C$ be a projective irreducible smooth curve, \m{n\geq 2} an integer and
\m{C_n} a primitive multiple curve of multiplicity $n$, with underlying smooth
curve $C$ and associated line bundle $L$ on $C$. Let $S$ be a smooth curve,
\m{P\in C} and \ \m{\pi:\kc\to S} \ a maximal reducible deformation of
\m{C_n}, with \ \m{\pi^{-1}(P)=C_n} (cf. \ref{Maxred_dr7}). Let
\m{\kc_1,\ldots,\kc_n} be the irreducible
components of $\kc$ and \ \m{d=-\deg(L)}. We suppose that \m{d>0}, i.e. $\pi$
is not a fragmented deformation, so for every \m{s\in S\backslash\{P\}} and
\m{i,j} distinct integers in \m{\{1,\ldots,n\}}, the two components
\m{\kc_{i,s}}, \m{\kc_{j,s}} of \m{\kc_s} intersect transversally in $d$
distinct points. As in \ref{Maxred2}, we denote also \m{t\circ\pi_i} by
\m{\pi_i}.

\sepprop

\begin{subsub}\label{notat}Intersections of components (notations) --  \rm
Let \m{\kz\subset\kc} be the closure in $\kc$ of the locus of the intersection
points of the components of \m{\kc_z=\pi^{-1}(z)}, \m{z\not=P}. Since $S$ is a
curve, $\kz$ is a curve. It intersects $C$ in a finite number of points.

For every proper subset \ \m{I\subset\{1,\ldots,n\}}, let \m{I^c} denote the
complement of $I$. If \m{J\subset\{1,\ldots,n\}} is such that
\m{I\subset J}, let \m{\kz^J_I} denote the subvariety of $\kz$ which is the
closure in $\kc$ of the locus of intersection points of the components
\m{\kc_{i,z}} and \m{\kc_{j,z}} of \m{\kc_z}, where \m{i\in I} and \m{j\in
J\backslash I}. Of course \m{\kz^J_I} is a union of components of \m{\kc_I\cap
\kc_{J\backslash I}} (the only other component being a curve with associated
reduced curve $C$). If \m{i,j\in\{1,\ldots,n\}} are distinct, we will note
more simply \ \m{\kz_{\{i\}}^{\{i,j\}}=\kz_{ij}}. Similarly we will note \
\m{\kz_I=\kz_I^{\{1,\ldots,n\}}}.

If \ \m{\emptyset\not=I\subsetneq J\subset\{1,\ldots,n\}}, let \m{\ki(I,J)}
denote the ideal sheaf of \m{\kz_I^J} in \m{\ko_{\kc_I}}. If
\m{J=\{1,\ldots,n\}}, we note \ \m{\ki(I)=\ki(I,\{1,\ldots,n\})}.
\end{subsub}

\sepprop

\begin{subsub}\label{spectrum} Spectrum -- \rm Let \m{i,j} be distinct integers
such that \m{1\leq i,j\leq n}, and \m{I=\{i,j\}}. According to \ref{Maxred2}
there exists a unique integer \m{p>0} such that \m{\ki_{C,\kc_I}/(\pi)} is
generated at any point \m{x\in C} by the image of an element of the form
\m{(\pi_i^p\alpha,0)}, where \m{\alpha\in\ko_{\kc_i,x}} is such that
\m{\alpha_{\mid C}\not=0} (and also by the image of an element
\m{(0,\pi_j^p\beta)}, where \m{\beta\in\ko_{\kc_j,x}} is such that
\m{\beta_{\mid C}\not=0}). Recall that $p$ is the smallest integer $q$ such
that \m{\ki_{C,\kc_I}} contains a non zero element of the form
\m{(\pi_i^q\lambda,0)} (or \m{(0,\pi_j^q\mu)}), with \m{\lambda_{\mid
C}\not=0} (resp. \m{\mu_{\mid C}\not=0}). Let
\[p_{ij} \ = \ p_{ji} \ = \ p ,\]
and \m{p_{ii}=0} for \m{1\leq i\leq n}. The symmetric matrix
\m{(p_{ij})_{1\leq i,j\leq n}} is called the {\em spectrum} of $\kc$.
\end{subsub}

\sepprop

\begin{subsub}\label{genecl} Generators of \m{(\ki_C^p+(\pi))/
(\ki_C^{p+1}+(\pi))} - \rm Let \m{i,j\in\{1,\ldots,n\}} be such that
\m{i\not=j} and \m{x\in C}. Since \m{\kc_{\{i,j\}}\subset\kc} there exists an
element \ \m{{\bf u}_i=(v_m)_{1\leq m\leq n}} \ of \m{\ko_{\kc,x}} such
that \m{v_i=0} and that the image of \m{(0,v_j)} generates \
\m{\ki_{C,\kc_{\{i,j\}},x}/(\ki_{C,\kc_{\{i,j\}},x}^2+(\pi))}. According to
\ref{p3}, the image of \m{{\bf u}_i} generates \m{\ki_C/(\ki_C^2+(\pi))} at
$x$, and for any \m{k\in\{1,\ldots,n\}} such that \m{k\not=i}, the image of
\m{(0,v_k)} generates \
\m{\ki_{C,\kc_{\{i,k\}},x}/(\ki_{C,\kc_{\{i,k\}},x}^2+(\pi))}. Hence from
theorem \ref{th1} there exists \ \m{\alpha_{ik}\in\ko_{\kc_k,x}} \ such that \
\m{v_k=\alpha_{ik}\pi_k^{p_{ik}}}, and
\[(\alpha_{ik}) \ = \ \ki_{\kz_{ik},\kc_k,x} \ . \]
Let $I$ be a proper subset of \m{\{1,\ldots,n\}}, $p$ its number of elements,
and
\[{\bf u}_I \ = \ \prod_{i\in I}{\bf u}_i \ . \]
Then the image of \m{{\bf u}_I} is a generator of \
\m{(\ki_{C,x}^p+(\pi))/(\ki_{C,x}^{p+1}+(\pi))}.
\end{subsub}

\sepprop

\begin{subsub}\label{prop2}{\bf Proposition: }{\bf 1 -- } The ideal sheaf \
\m{\ki_{\kc_I,\kc}} \ of \m{\kc_I} in $\kc$ at $x$ is generated by \m{{\bf
u}_I}. In particular the ideal sheaf of \m{\kc_i} in $\kc$ at $x$ is generated
by \m{{\bf u}_i}.

{\bf 2 -- } \m{\ki_{\kc_I,\kc}} is a line bundle on \m{\kc_{I^c}}.
\end{subsub}

Of course {\bf 2-} is a consequence of {\bf 1-}. The proof of {\bf 1-} is
similar to that of proposition 4.3.3 in \cite{dr7}.

In particular, for every \m{i\in\{1,\ldots,n\}}, let \
\m{J_i=\{1,\ldots,n\}\backslash\{i\}}. Then the ideal sheaf of \
\m{\kc_{J_i}\subset\kc} \ at $x$ is generated by
\[{\bf u}_{J_i} \ = \ (0,\ldots,0,A_i\pi_i^{q_i},0,\ldots,0) \ , \]
with
\[A_i \ = \ \prod_{1\leq k\leq n,k\not=i}\alpha_{ki}, \qquad q_i \ = \
\sigg_{1\leq k\leq n}p_{ik} \ . \]

\sepprop

\begin{subsub}\label{prop2b}{\bf Proposition: } Let \ \m{{\bf
u}_I=(w_i)_{1\leq i\leq n}}. For every \ \m{j\in\{1,\ldots,n\}\backslash
I}, let \m{w_j=v_j\pi_j^{m_j}}, with \m{v_{j|C}\not=0}. Then \m{v_j} is a
generator of the ideal of \m{\kz_{\{j\}}^{I\cup\{j\}}} at $x$.
\end{subsub}

\begin{proof} This follows from corollary \ref{coro1}, 2-.
\end{proof}

\sepprop

\begin{subsub}\label{prop3}{\bf Proposition: } {\bf 1 -- } Let \m{i,j,k} be
distinct integers such that \m{1\leq i,j,k\leq n}. Then if \m{p_{ij}<p_{jk}},
we have \m{p_{ik}=p_{ij}}.

{\bf 2 -- } Let \m{i\in\{1,\ldots,n\}}. Then we have \
\m{\ki_{C,x}=({\bf u}_i)+(\pi)}.

{\bf 3 -- } Let \m{i\in\{1,\ldots,n\}} and \m{v=(v_m)_{1\leq m\leq
n}\in\ki_{C,x}} such that \m{v_i} is a multiple of \m{\pi_i^p}, with \m{p>0}.
Then we have \ \m{v\in({\bf u}_i)+(\pi^p)}
\end{subsub}

The proof is similar to that of propositions 4.3.4 and 4.3.5 in \cite{dr7}.

\sepprop

\begin{subsub}\label{const1}{Construction by induction on $n$ -- } \rm
Let $i$ be an integer such that \m{1\leq i\leq n}.
Let \m{\kb_i} be the image of \m{\ko_{\kc_{J_i}}} in the sheaf of
$\C$-algebras on $C$ which at any point $x$ is \ \m{\prod_{1\leq j\leq
n,j\not=i}\ko_{\kc_j}/(A_j\pi_j^{q_j})}; it is also a sheaf of $\C$-algebras
on $C$. For every point $x$ of $C$ and every \m{\alpha=(\alpha_m)_{1\leq m\leq
n}} in \ \m{\prod_{1\leq j\leq n}\ko_{\kc_j,x}}, we denote by
\m{b_{i,x}(\alpha)} its image in \ \m{\prod_{1\leq j\leq
n,j\not=i}\ko_{\kc_j,x}/(A_j\pi_j^{q_j})} (obtained by forgetting the $i$-th
coordinate of $\alpha$).

We denote by \m{\ka_i} the sheaf of $\C$-algebras on $C$ which at any point
$x$ is \ \m{\ko_{\kc_i}/(A_i\pi_i^{q_i})}.

If \m{1\leq k\leq n}, \m{x\in C} and \m{\alpha\in\ko_{\kc_k,x}}, let
\m{[\alpha]} denote the image of $\alpha$ in
\m{\ko_{\kc_k,x}/(A_k\pi_k^{q_k})}.
\end{subsub}

\sepprop

\begin{subsub}\label{ecl9}{\bf Proposition: } There exists a morphism of
sheaves of $\C$-algebras on $C$
\[\Phi_i : \kb_i\lra\ka_i\]
such that for every point $x$ of $C$ and all \m{(\alpha_m)_{1\leq m\leq
n,m\not=i}\in\ko_{\kc_{J_i},x}}, \m{\alpha_i\in\ko_{\kc_i,x}}, we have \Nligne
\m{\alpha=(\alpha_m)_{1\leq m\leq n}\in\ko_{\kc,x}} if and only if \
\m{\Phi_{i,x}(b_{i,x}(\alpha))=[\alpha_i]}.
\end{subsub}

The proof is similar to that of proposition 4.4.1 of \cite{dr7}.

\sepprop

For \m{1<i\leq n}, let
\[A^{[i]} \ = \ \prod_{1\leq k< i}\alpha_{ki} \ . \]
The following result is a generalization of corollary 4.4.3 of \cite{dr7}:

\sepprop

\begin{subsub}\label{coro2}{\bf Corollary: } Let $N$ be an integer such that
\m{N\geq\max_{1\leq i\leq n}(q_i)}. Let \m{x\in C},\Nligne
\m{\beta\in\ko_{\kc_1,x}\times\cdots\ko_{\kc_n,x}} \ and \m{u=(u_i)_{1\leq
i\leq n}\in\ko_{\kc,x}} such that \m{u_i} and \m{A^{[i]}} are relatively prime
in \m{\ko_{\kc_i,x}}, for \m{1< i\leq n}. Suppose that \ \m{\lbrack\beta
u\rbrack\in\ko_{\kc,x}/(\pi^N)}. Then we have \
\m{\lbrack\beta\rbrack\in\ko_{\kc,x}/(\pi^N)}.
\end{subsub}

\begin{proof} By induction on $n$. It is obvious if $n=1$. Suppose that the
lemma is true for \m{n-1}. Let \ \m{I=\{1,\ldots,n-1\}}. So we have \
\m{[\beta_{\mid\kc_1\times\cdots\kc_{n-1}}]
\in\ko_{\kc_I,x}/(\pi_1,\ldots,\pi_{n-1})^N} by the induction hypothesis.
Let $\gamma$ (resp. $v$) be the image of \m{\beta} (resp. $u$) in \m{\kb_n}. To
show that \m{\lbrack\beta\rbrack\in\ko_{\kc,x}/(\pi^N)} it is enough to
verify that
\[\Phi_n(\gamma) \ = \ \lbrack\beta_n\rbrack .\]
We have \ \m{\Phi_n(\gamma v) \ = \ \lbrack\beta_nu_n\rbrack} because
\m{\lbrack\beta u\rbrack\in\ko_{\kc,x}/(\pi^N)}, and \m{\Phi_n(v)=\lbrack u_n
\rbrack} because \ \m{u\in\ko_{\kc,x}}. So we have
\[\Phi_n(\gamma)\lbrack u_n\rbrack \ = \ \Phi_n(\gamma)\Phi_n(v) \ = \
\Phi_n(\gamma v) \ = \ \lbrack\beta_nu_n\rbrack \ = \
\lbrack\beta_n\rbrack\lbrack u_n\rbrack .\]
Since \m{u_{\mid C}\not=0}, \m{\lbrack u_n\rbrack} is not a zero divisor
in \m{\ko_{\kc_n,x}/(A_n\pi_n^{q_n})}, so we have \ \m{\Phi_n(\gamma)=
\lbrack\beta_n\rbrack}.
\end{proof}

\sepprop

The following result is a generalization of proposition 4.4.7 of \cite{dr7}:

\sepprop

\begin{subsub}\label{ecl10}{\bf Proposition: } Let \ \m{(\alpha_1\pi_1^{m_1},
\ldots,\alpha_n\pi_n^{m_n})\in\ko_{\kc,x}} \ be such that for \m{1<i\leq n},
\m{\alpha_{i|C}\not=0} and that \m{\alpha_i} and \m{A^{[i]}} are relatively
prime. Let \m{\beta=(\beta_1,\ldots,\beta_n)\in\ko_{\kc,x}} be such that for
\m{1\leq i\leq n}, \m{\beta_i} is a multiple of \m{\alpha_i} in
\m{\ko_{\kc_i,x}}. Let \ \m{M=m_1+\cdots m_n}. Then
\[\left(\frac{\beta_1}{\alpha_1}\pi_1^{M-m_1},\ldots,
\frac{\beta_n}{\alpha_n}\pi_n^{M-m_n}\right) \ \in \ \ko_{\kc,x} \ .\]
\end{subsub}

\begin{proof} By induction on $n$. It is obvious for \m{n=1}. Suppose that it
is true for \m{n-1\geq 1}. Let \ \m{I=\{1,\ldots,n-1\}}. Then \
\m{(\alpha_1\pi_1^{m_1},\ldots,\alpha_{n-1}\pi_{n-1}^{m_{n-1}})\in
\ko_{\kc_{I,x}}}. Hence, by the induction hypothesis, we have
\[\big(\frac{\beta_1}{\alpha_1}\pi_1^{M-m_1-m_n},\ldots,
\frac{\beta_{n-1}}{\alpha_{n-1}}\pi_{n-1}^{M-m_{n-1}-m_n}
\big) \ \in \ \ko_{\kc_{I,x}} \ .\]
So there exists \m{\gamma\in\ko_{\kc_{n,x}}} such that
\[u \ = \ \big(\frac{\beta_1}{\alpha_1}\pi_1^{M-m_1-m_n},\ldots,
\frac{\beta_{n-1}}{\alpha_{n-1}}\pi_{n-1}^{M-m_{n-1}-m_n},\gamma\big)
\ \in \ \ko_{\kc,x} \ .\]
Multiplying by  \m{(\alpha_1\pi_1^{m_1},\ldots,\alpha_n\pi_n^{m_n})} we see
that \ \m{(\beta_1\pi_1^{M-m_n},\ldots,\beta_{n-1}\pi_{n-1}^{M-m_n},\gamma
\alpha_n\pi_n^{m_n})\in\ko_{\kc,x}}. Subtracting \m{\beta\pi^{M-m_n}}, we find
that \
\m{(0,\ldots,0,\gamma\alpha_n\pi_n^{m_n}-\beta_n\pi_n^{M-m_n})\in\ko_{\kc,x}}.
By proposition \ref{prop2}, the ideal of \m{\kc_I\subset\kc} at $x$ is
generated by \m{(0,\ldots,0,A_n\pi^{q_n})}. Hence we can write
\[\gamma\alpha_n\pi_n^{m_n}-\beta_n\pi_n^{M-m_n} \ = \ \lambda A_n\pi_n^{q_n}\]
for some \m{\lambda\in\ko_{\kc_n,x}}. Since \m{\alpha_n} divides \m{\beta_n},
it divides also \m{\lambda A_n}, and since \m{\alpha_n} and \m{A_n} are
relatively prime, \m{\alpha_n} divides $\lambda$, and we have \
\m{\displaystyle\gamma\pi_n^{m_n}-\frac{\beta_n}{\alpha_n}\pi_n^{M-m_n}=\mu
A_n} \ for some \m{\mu\in\ko_{\kc_{n,x}}}. We have \
\m{v=(0,\ldots,0,\mu A_n\pi_n^{q_n})\in\ko_{\kc,x}}.
Now we have
\[\pi^{m_n}u-v \ = \
\left(\frac{\beta_1}{\alpha_1}\pi_1^{M-m_1},\ldots,\frac{\beta_n}{\alpha_n}
\pi_n^{M-m_n}\right) \ \in \ \ko_{\kc,x} \ . \]
\end{proof}

\sepprop

\begin{subsub}\label{frag}Reducible deformations and fragmented deformations --
\rm We show here that it is not true that there always exists a fragmented
deformation above a reducible one (one could think that a reducible deformation
could be obtained simply by gluing some curves together in a fragmented one).

Let $C$ be a smooth curve of degree \m{d>0} in \m{\P_2}, and \m{f\in
H^0(\ko_{\P_2}(d))} an equation of $C$. Let \m{h\in H^0(\ko_{\P_2}(d))},
\[\kc(h) \ = \ \{(x,t)\in\P_2\times\C;(f+th)(x)=0\} \]
and \ \m{\pi_h:\kc(h)\to\C} \ the projection. Now let \m{h_1,\ldots,h_n\in
H^0(\ko_{\P_2}(d))}. Let \Nligne
\m{\kd=\kc(h_1)\cup\cdots\cup\kc(h_n)\subset\P_2\times\C} \ and \
\m{\pi=(\pi_1,\ldots,\pi_n):\kd\to\C}. If \m{h_1,\ldots,h_n} are sufficiently
general, there exists a suitable neighborhood $U$ of $0$ in $\C$ such that
\m{\kc=\pi^{-1}(U)} is a reducible deformation of $C$. Of course
\m{\pi^{-1}(0)} is the multiple curve in \m{\P_2} defined by \m{f^n}. Now we
have \m{p_{ij}=1} for \m{1\leq i<j\leq n}, and
\[{\bf u}_1 \ = \ (0,(h_2-h_1)t,\ldots,(h_n-h_1)t) \ . \]
Let \m{x\in C} and suppose that \m{h_2-h_1} vanishes at $x$, but not
\m{h_3-h_1}. We will show that there is no fragmented deformation above $\kc$,
i.e. there is no fragmented deformation \m{\kc'} of $C$ with components
\m{\kc(h_1)_{|U}},$\ldots$, \m{\kc(h_n)_{|U}} such that \
\m{\ko_\kc\subset\ko_{\kc'}}. If such a fragmented deformation exists, its
spectrum is the same as that of $\kc$, and the ideal of \m{\kc(h_1)_{|U}} in
\m{\kc'} is generated at $x$ by an element
\[{\bf u}' \ = \ (0,t,a_2t,\ldots,a_nt)\]
with \m{a_i} invertible in \m{\ko_{\kc(h_i),x}} for \m{2\leq i\leq n}. Then
\m{{\bf u}'} is a multiple of $\bf u$, i.e there exists
\m{\alpha\in\ko_{\kc,x}} of the form
\[\alpha \ = \ (\alpha_1,h_2-h_1,\frac{1}{a_3}(h_3-h_1),\ldots,
\frac{1}{a_n}(h_n-h_1)) \ . \]
But this is impossible because \ \m{\dsp(h_2-h_1)_{|C}\not=
\frac{1}{a_3}(h_3-h_1)_{|C}}.
\end{subsub}

\sepsub

\end{sub}

\sepsub

\Ssect{Localization and construction of maximal reducible deformations}{argn}

We keep the notations of \ref{Mnprop}. The following result generalizes
proposition \ref{prop1b}:

\sepprop

\begin{subsub}\label{ecl11}{\bf Proposition: } Let \
\m{u=(\alpha_1\pi_1^{m_1},\ldots,\alpha_n\pi_n^{m_n})\in\ko_{\kc,x}} , with
\m{m_i>q_i} for \m{1\leq i\leq n}. Then $u$ is a multiple of $\pi$ in
\m{\ko_{\kc,x}} , i.e. \
\m{(\alpha_1\pi_1^{m_1-1},\ldots,\alpha_n\pi_n^{m_n-1})\in\ko_{\kc,x}}.
\end{subsub}

\begin{proof} By induction on $n$. The result is true for \m{n=1}.
Suppose that \m{n>1} and that it is true for \m{n-1}. Let \
\m{I=\{0,\ldots,n-1\}}. Then \
\m{u'=(\alpha_1\pi_1^{m_1-1},\ldots,\alpha_{n-1}\pi_{n-1}^{m_{n-1}-1})
\in\ko_{\kc_I,x}} (by the induction hypothesis). Let \m{v\in\ko_{C,x}} be such
that its image in \m{\ko_{\kc_I,x}} is \m{u'}, i.e \m{v=(v_i)_{1\leq i\leq n}}
and \ \m{v_i=\alpha_i\pi_1^{m_i-1}} \ for \m{1\leq i\leq n-1}. Then \ \m{u-\pi
v\in\ki_{\kc_I,\kc,x}}, hence there exists \m{\lambda\in\ko_{\kc_n,x}} such that
\[u-\pi v \ = \ (0,\ldots,0,\lambda A_n\pi_n^{q_n}) \ . \]
We have only to show that $\lambda$ is divisible by \m{\pi_n}. We have \
\m{v_n=\mu\pi_n^{q_n-1}}, with \Nligne \m{\mu=\pi_n^{m_n-q_n}a_n-\lambda A_n},
so we have to show that $\mu$ is divisible by \m{\pi_n}. Let
\m{\tau=(\tau_i)_{1\leq i\leq n}\in\ki_{\kz,x}}, such that \m{\tau_{|C}\not=0}.
Then we have
\[\tau v \ = \ (\pi_1^{m_1-1}a_1\tau_1,\ldots,
\pi_{n-1}^{m_{n-1}-1}a_{n-1}\tau_{n-1},\pi_n^{q_n-1}\mu\tau_n) \ . \]
Since \m{m_i-1\geq q_i} and \m{A_i} divides \m{\tau_i} for \m{1\leq i\leq n-1},
we have
\[(\pi_1^{m_1-1}a_1\tau_1,\ldots,\pi_{n-1}^{m_{n-1}-1}a_{n-1}\tau_{n-1},0)\in
\ko_{\kc,x} \ . \]
Hence \ \m{w=(0,\ldots,0,\pi_n^{q_n-1}\mu\tau_n)\in\ko_{\kc,x}}. Hence $w$ is a
multiple of \m{(0,\ldots,0,\pi_n^{q_n}A_n)}. Since \m{\tau_n} is not divisible
by \m{\pi_n}, $\mu$ is divisible by \m{\pi_n}.
\end{proof}

\sepprop

Let \m{I=\{1,\ldots,n-1\}} and \m{x\in C}. Let \m{\boldsymbol{B_n}} be the
image of \m{\ko_{\kc_{I,x}}} in \Nligne
\m{\prod_{i=1}^{n-1}\ko_{\kc_i,x}/\big((A_i\pi_i^{q_i})+(\pi_i^{q_i+1})\big)
\times\ldots\times\ko_{\kc_{n-1},x}/\big((A_{n-1}\pi_{n-1}^{q_{n-1}})+
(\pi_{n-1}^{q_{n-1}+1})\big)}. As in \ref{local}, we can deduce from the
morphism \m{\Phi_{n-1}} of proposition \ref{ecl9} the morphism
\[\ov{\Phi}_{n-1}:\boldsymbol{B_n}\lra\ko_{\kc_n,x}/\big((A_n\pi_n^{q_n})+
(\pi_n^{q_n+1})\big) \ . \]
We can then define a new sheaf of algebras $\boldsymbol{\ka}$ on $C$ by
\[\boldsymbol{\ka}_x \ = \ \{(u,v)\in\boldsymbol{B_n}\times
\ko_{\kc_n,x}/\big((A_n\pi_n^{q_n})+(\pi_n^{q_n+1})\big) \ ; \
\ov{\Phi}_{n-1}(u)=v\} \ . \]
The proof of the following result is similar to that of proposition
\ref{prop1c}:

\sepprop

\begin{subsub}\label{prop1cn}{\bf Proposition: } The natural projection \
\m{\ko_{\kc,x}\to\ka_x} \ induces an isomorphism
\[\ko_{\kc,x}/(\pi) \ \simeq \ \boldsymbol{\ka}_x/(\pi) \ . \]
\end{subsub}

It is then possible to define as for double primitive curves in
\ref{loc_red_def} the notion of a {\em local reducible deformation} of a
primitive multiple curve of multiplicity $n$. By proposition \ref{prop1cn}, a
global reducible deformation induces a local one. And as for double curves, if
a primitive multiple curve of multiplicity $n$ has a local reducible
deformation, it has also a global one.

\end{sub}

\end{document}